\documentstyle{article}

\newcommand{\bigzerou}{%
\smash{\lower1.7ex\hbox{\bg 0}}}
\newtheorem{prop}{Proposition}
\newtheorem{defi}{Definition}

\newtheorem{conj}{Conjecture}

\newtheorem{Rem}{Remark}
\newtheorem{Q}{Question}

\newcommand{\ba}{\begin{eqnarray}}
\newcommand{\ea}{\end{eqnarray}}
\newcommand{\no}{\nonumber}

\newcommand{\mapright}[1]{%
\smash{\mathop{%
\hbox to 1.0cm{\rightarrowfill}}\limits^{#1}}}
\newcommand{\mapleft}[1]{%
\smash{\mathop{%
\hbox to 1.3cm{\leftarrowfill}}\limits^{#1}}}

\begin{document}
\title{
\begin{flushright}
  \begin{minipage}[b]{5em}
    \normalsize
    ${}$      \\
  \end{minipage}
\end{flushright}
{\bf Virtual Gromov-Witten Invariants and the Quantum Cohomology Rings of
General Type Projective Hypersurfaces}}
\author{Masao Jinzenji\\
\\ 
\it Graduate School of Mathematical Sciences\\
\it University of Tokyo\\
\it  Meguro-ku, Tokyo 153-8914, Japan\\}
\maketitle
\begin{abstract}
In this paper, we propose another characterization of the generalized
mirror transformation on the quantum cohomology rings of general type 
projective hypersurfaces. This characterics is useful for explicit 
determination of the form of the generalized mirror transformation. 
As applications, we rederive the generalized mirror transformation up to 
$d=3$ rational Gromov-Witten invariants obtained in our previous
article, and determine explicitly the the generalized mirror 
transformation for the $d=4, 5$ rational Gromov-Witten invariants in the 
case when the first Chern class of the hypersurface equals $-H$ (i.e.,
$k-N=1$).
\footnote{ e-mail 
address: jin@hep-th.phys.s.u-tokyo.ac.jp}  
\end{abstract}
\section{Introduction and Statement of the Main Results}
Recently, some works on the quantum cohomology ring of the general type
projective hypersurface have appeared \cite{yau}, \cite{gath},
\cite{gene}.
In \cite{yau}, Lian, Liu and Yau generalized their mirror principle 
to the case of the general type projective hypersurface and proposed 
a theoretical recipe to construct the generating function of a certain 
type of Gromov-Witten invariants including gravitational descendants 
from the hypergeometric data.
In \cite{gath}, Gathmann considered the relative Gromov-Witten
invariants of the projective space with various tangency conditions on 
the ample hypersurface in the projective space. He verified the
recursive formulas which increase the tangency condition like the
results of Caporaso and Harris  \cite{harris}, and proposed an algorithm 
of computing the Gromov-Witten invariants of the hypersurface as a limit 
of iterative application of the recursive formulas. 
    
In this paper, we continue the analysis of \cite{gene} on 
the quantum K\"ahler sub-ring $QH^{*}_{e}(M_{N}^{k})$, 
where $M_{N}^{k}$ is the degree $k$ hypersurface 
in $CP^{N-1}$, especially in the
case when its first Chern class is negative.  
Our approach is different from the ones of \cite{yau} and \cite{gath}, 
mainly because we don't use Mumford-Morita class (gravitational
descendants). In \cite{gene}, we 
proposed the generalized mirror transformation on the quantum cohomology 
rings of $M_{N}^{k},\;\;(N-k\leq 1)$ that represents the sturctural
constant $L_{n}^{N,k,d}$ of $QH^{*}_{e}(M_{N}^{k})$ in terms of 
the virtual sturctural constants $\tilde{L}_{n}^{N,k,d}$, which is the
analogue of the hypergeometric series used in the mirror calculation in
our context.

Now we restate the main conjecture in \cite{gene}.  
Let $P_{m}$ be the set of partitions of $m$ into positive 
integers and $\sigma_{m}$ be an element of $P_{m}$. We also denote 
the length of a partition $\sigma_{m}$ by $l(\sigma_{m})$ (i.e., 
$\sigma_{m}:m=d_{1}+d_{2}+
\cdots+d_{l(\sigma_{m})},\;\;d_{1}\geq d_{2}\geq\cdots
\geq d_{l(\sigma_{m})}\geq1$). We denote by $\mbox{mul}(i,\sigma_{m}),\;
(1\leq i \leq m)$ the multiplicity of $i$ in $\sigma_{m}$. 
Our previous conjecture is the following:
\begin{conj}
The generalized mirror transformation takes the form
\begin{eqnarray}
L^{N,k,d}_{n}&=&\sum_{m=0}^{d-1}\sum_{\sigma_{m}\in P_{m}}
(-1)^{l(\sigma_{m})}\frac{d^{l(\sigma_{m})}}
{\prod_{j=1}^{l(\sigma_{m})}d_{j} \prod_{i=1}^{m}\mbox{mul}(i,\sigma_{m})!}
\cdot\prod_{i=1}^{l(\sigma_{m})}
\tilde{L}_{1+(k-N)d_{i}}^{N,k,d_{i}}\cdot G_{d-m}^{N,k,d}(n;\sigma_{m}),\no\\
\label{gene}
\end{eqnarray}
where $G_{d-m}^{N,k,d}(n;\sigma_{m})$ is a polynomial of 
$\tilde{L}_{n}^{N,k,d}$  with weighted degree $d$. 
\end{conj}
Now, we propose a conjecture on the explicit form of 
$G_{d-m}^{N,k,d}(n;\sigma_{m})$.  
\begin{conj} 
\begin{eqnarray}
V_{d-m}^{N,k,d}(n;\sigma_{m})&:=& 
\frac{1}{k}v({\cal O}_{e^{N-2-n}}{\cal O}_{e^{n-1-(k-N)d}}
\prod_{i=1}^{l(\sigma_{m})}{\cal O}_{e^{1+(k-N)d_{i}}})_{d-m}\cdot
\frac{1}{(d-m)^{l(\sigma_{m})-1}}\no\\
&\simeq& G_{d-m}^{N,k,d}(n;\sigma_{m})
\label{main}
\end{eqnarray}
where $v({\cal O}_{e^{N-2-n}}{\cal O}_{e^{n-1-(k-N)d}}
\prod_{i=1}^{l(\sigma_{m})}{\cal O}_{e^{1+(k-N)d_{i}}})_{d-m}$
is the virtual Gromov-Witten invariant defined below. If
 $\l(\sigma_{m})\leq 1$ or $d-m=1$, (\ref{main}) becomes an exact equality. 
\end{conj}
We have to make some remarks on the meaning of $\simeq$ in (\ref{main}).
In the case when $l(\sigma_{m})\leq 1$ or $d-m=1$, 
we can see the coincidence between 
$G_{d-m}^{N,k,d}(n;\sigma_{m})$ and $V_{d-m}^{N,k,d}(n;\sigma_{m})$, but 
when $l(\sigma_{m})\geq 2$ and $d-m\geq
2$, we have to modify slightly the
coefficients of some polynomials in $\tilde{L}_{n}^{N,k,d}$ that appear in
$V_{d-m}^{N,k,d}(n;\sigma_{m})$. These cases indeed appear when $d\geq 4$. 
Even in such 
situations, $V_{d-m}^{N,k,d}(n;\sigma_{m})$ strongly predicts the form
of $G_{d-m}^{N,k,d}(n;\sigma_{m})$.  

Now we turn into the definition of the virtual Gromov-Witten invariants, 
that are the main ingredients of this paper. 
 
\begin{defi}
The virtual Gromov-Witten invariant
$v(\prod_{j=1}^{n}{\cal O}_{e^{a_{j}}})_{d}$ on $M_{N}^{k}$ is the 
rational number that satisfy the condition:\\
(i) initial condition
\begin{eqnarray}
&&v({\cal O}_{e^{a}}{\cal O}_{e^{b}}{\cal O}_{e^{c}})_{0}=
k\cdot\delta_{a+b+c,N-2},\no\\
&&v(\prod_{j=1}^{n}{\cal O}_{e^{a_{j}}})_{0}=0,
\;\;(n\neq 3),\no\\
&&\frac{1}{k}v({\cal O}_{e^{N-2-n}}{\cal O}_{e^{n-1-(k-N)d}}{\cal O}_{e})_{d}=
\tilde{L}_{n}^{N,k,d}-\tilde{L}_{1+(k-N)d}^{N,k,d},\;\;(d\geq 1),
\label{initial}
\end{eqnarray}
(ii) flat metric condition
\begin{eqnarray}
&&v({\cal O}_{e^{0}}{\cal O}_{e^{a}}{\cal O}_{e^{b}})_{0}
=k\cdot\delta_{a+b,N-2},\no\\
&&v({\cal O}_{e^{0}}\prod_{j=1}^{n}{\cal O}_{e^{a_{j}}})_{d}=0,
\;\;(d\geq 1,\;\;\mbox{or} \;\;\;d=0,\;\; n\neq 2),
\end{eqnarray}
(iii) topological selection rule\\
\begin{equation}
v(\prod_{j=1}^{n}{\cal O}_{e^{a_{j}}})_{d}\neq0
\Longrightarrow (N-5)+(N-k)d=\sum_{j=1}^{n}(a_{j}-1),
\end{equation}
(iv) K\"ahler equation\\
\begin{equation}
v({\cal O}_{e}\prod_{j=1}^{n}{\cal O}_{e^{a_{j}}})_{d}=
d\cdot v(\prod_{j=1}^{n}{\cal O}_{e^{a_{j}}})_{d},
\end{equation}
(v) associativity equation \cite{km}, \cite{rt}\\
\begin{eqnarray}
&&\sum_{d_{1}=0}^{d}\sum_{\{\alpha_{*}\}\coprod\{\beta_{*}\}=\{n_{*}\}}\sum_{i=0}^{N-2}v({\cal O}_{e^{a}}{\cal O}_{e^{b}}
(\prod_{\alpha_{j}\in\{\alpha_{*}\}}{\cal O}_{e^{\alpha_{j}}}){\cal O}_{e^{i}})_{d_{1}}
v({\cal O}_{e^{N-2-i}}  
(\prod_{\beta_{j}\in\{\beta_{*}\}}{\cal O}_{e^{\beta_{j}}}){\cal O}_{e^{c}}
{\cal O}_{e^{d}})_{d-d_{1}}\no\\
&&=\sum_{d_{1}=0}^{d}\sum_{\{\alpha_{*}\}\coprod\{\beta_{*}\}=\{n_{*}\}}\sum_{i=0}^{N-2}
v({\cal O}_{e^{a}}{\cal O}_{e^{c}}
(\prod_{\alpha_{j}\in\{\alpha_{*}\}}{\cal O}_{e^{\alpha_{j}}}){\cal O}_{e^{i}})_{d_{1}}
v({\cal O}_{e^{N-2-i}}  
(\prod_{\beta_{j}\in\{\beta_{*}\}}{\cal O}_{e^{\beta_{j}}}){\cal O}_{e^{b}}
{\cal O}_{e^{d}})_{d-d_{1}},\no\\
&&(a+b+c+d+\sum_{j=1}^{m}(n_{j}-1)=N-2+(N-k)d).
\end{eqnarray}
\end{defi}
The difference between the virtual Gromov-Witten invariants and the 
ordinary Gromov-Witten invariants comes from the initial condition
(\ref{initial}). In ordinary cases, we set
\begin{equation}
\frac{1}{k}\langle
{\cal O}_{e^{N-2-n}}{\cal O}_{e^{n-1-(k-N)d}}{\cal O}_{e}\rangle_{d}=
{L}_{n}^{N,k,d},\;\;(d\geq 1).
\label{ord}
\end{equation}
(\ref{initial}) and (\ref{ord}) indeed differ when $d>1$.
With the above conditions, we can completely determine the virtual 
Gromov-Witten invariants like the ordinary Gromov-Witten invariants of
genus $0$.  
Using the generating function of the virtual Gromov-Witten invariants:
\begin{eqnarray}
F_{v}(z,t_{i})=\frac{1}{3!}\sum_{i,j,k=0}^{N-2}
k\cdot\delta_{i+j+k,N-2}\cdot t_{i}t_{j}t_{k}+f_{v}(z,t_{i}),
\end{eqnarray}
we can state the main conjecture in a more compact form,
\begin{equation}
L_{n}^{N,k,d}\simeq\frac{1}{k}\int_{C_{0}}dz\frac{1}{z^{d+1}}
\exp(-d\sum_{j=1}^{\infty}\frac{\tilde{L}_{1+(k-N)j}^{N,k,j}}{j}
\frac{\partial_{1+(k-N)j}}{\partial_{1}}z^{j})\cdot \partial_{N-2-n}
\partial_{n-1-d(k-N)}\partial_{1}f_{v}(z,t_{i})|_{t_{i}=0},
\label{compact}
\end{equation}
where $t_{j}\; (j=0,\cdots N-2)$ is the variable that couples to 
the element ${\cal O}_{e^{j}}$ of $QH_{e}^{*}(M_{N}^{k})$ and $z$ is the 
formal degree counting variable.

This paper is organized as follows. In Section 2, we introduce the
notation of the quantum K\"ahler subring of $M_{N}^{k}$, and review the
results obtained in \cite{cj}, \cite{6} and \cite{gene}. In Section 3, 
we reproduce the formulas obtained in \cite{gene} under the assumption 
of Conjecture 2. In Section 4, we construct the explicit form of the 
generalized mirror transformation for degree $4, 5$ rational Gromov
-Witten invariants of $M_{k-1}^{k}$ using the Conjecture 2 and some
numerical data obtained from the fixed point computation in \cite{torus}.
\section{Quantum K\"ahler Sub-ring of Projective Hypersurfaces}
\subsection{Notation} 
In this section, we introduce the quantum K\"ahler sub-ring 
of the quantum cohomology ring of a degree $k$ hypersurface in
$CP^{N-1}$.
Let $M_{N}^{k}$ be a hypersurface of degree $k$ in $CP^{N-1}$.
 We denote by $QH^{*}_{e}(M_{N}^{k})$ the 
subring of the quantum cohomology ring $QH^{*}(M_{N}^{k})$
generated by ${\cal O}_{e}$ induced from the K\"ahler form $e$ 
(or, equivalently the intersection $H\cap M_{N}^{k}$ between a hyperplane
class $H$ of $CP^{N-1}$ and $M_{N}^{k}$).
Additive basis of $QH_{e}^{*}(M_{N}^{k})$ is given by 
${\cal O}_{e^{j}}\;\;(j=0,1,\cdots,N-2)$, which is induced from 
$e^{j}\in H^{j,j}(M_{N}^{k})$.
 The multiplication rule of $QH^{*}_{e}(M_{N}^{k})$ 
is determined by the Gromov-Witten invariant of genus $0$ 
$\langle {\cal O}_{e}{\cal O}_{{e}^{N-2-m}}
{\cal O}_{{e}^{m-1-(k-N)d}}\rangle_{d,M_{N}^{k}}$ and
it is given as follows:
\begin{eqnarray}
 L_{m}^{N,k,d} &:=&\frac{1}{k}\langle {\cal O}_{e}{\cal O}_{{e}^{N-2-m}}
{\cal O}_{{e}^{m-1-(k-N)d}}\rangle_{d},\no\\
\no\\
{\cal O}_{e}\cdot 1&=&{\cal O}_{e},\nonumber\\
{\cal O}_{e}\cdot{\cal O}_{{e}^{N-2-m}}&=&{\cal O}_{{e}^{N-1-m}}+
\sum_{d=1}^{\infty}L_{m}^{N,k,d}q^{d}{\cal O}_{{e}^{N-1-m+(k-N)d}},\no\\
q&:=&\exp(t),
\label{gm}
\end{eqnarray}
where the subscript $d$ counts the degree of the rational curves
measured by $e$. Therefore,  $q=\exp(t)$ is the degree counting 
parameter. 
\begin{defi}
We call $L_{n}^{N,k,d}$ the structural constant of weighted degree $d$.
\end{defi}
Since $M_{N}^{k}$ is a complex $(N-2)$ dimensional manifold, we see that
a structure constant $L_{m}^{N,k,d}$
is non-zero only if the following condition is satisfied:
\begin{eqnarray}
&& 1\leq N-2-m\leq N-2, 1\leq m-1+(N-k)d\leq N-2,\no\\
&\Longleftrightarrow &max\{0,2-(N-k)d\}\leq m \leq min\{N-3,N-1-(N-k)d\}.
\label{sel}
\end{eqnarray}
We rewrite (\ref{sel}) into 
\begin{eqnarray}
(N-k\geq 2) &\Longrightarrow& 0\leq m \leq (N-1)-(N-k)d\no\\
(N-k=1,d=1)&\Longrightarrow& 1\leq m \leq N-3\no\\
(N-k=1,d\geq2)&\Longrightarrow& 0\leq m \leq N-1-(N-k)d\no\\
(N-k\leq 0)&\Longrightarrow& 2+(k-N)d\leq m \leq N-3.
\label{flasel}
\end{eqnarray}
From (\ref{flasel}), we easily see that the number of the non-zero
structure constants $L_{m}^{N,k,d}$ is finite except for the case of $N=k$.
Moreover, if $N\geq 2k$, the non-zero structure constants come only from
the $d=1$ part and the non-vanishing $L_{m}^{N,k,1}$  
is determined by $k$ and  
independent of $N$. 
The $N\geq 2k$ region is studied by Beauville \cite{beauville}, 
and his result plays 
the role of an initial condition of our discussion later.
Explicitly, they are given by the formula :
\begin{equation}
\sum_{n=0}^{k-1}L_{n}^{N,k,1}w^{n}=k\prod_{j=1}^{k-1}(jw+(k-j)),
\label{one}
\end{equation}
and the other $L_{n}^{N,k,d}$'s all vanishes.
In the case of $N=k$, the multiplication rule of $QH^{*}_{e}(M_{k}^{k})$ is
given as follows:
\begin{eqnarray}
{\cal O}_{e}\cdot 1&=&{\cal O}_{e},\nonumber\\
{\cal O}_{e}\cdot{\cal O}_{{e}^{k-2-m}}&=&
(1+\sum_{d=1}^{\infty}q^{d}L_{m}^{k,k,d}){\cal O}_{{e}^{k-1-m}}
\;\;(m=2,3,\cdots,k-3),\no\\
{\cal O}_{e}\cdot{\cal O}_{{e}^{k-3}}  &=&{\cal O}_{e^{k-2}}.
\label{calabi}
\end{eqnarray}
We introduce here the generating function of the structure constants 
of the Calabi-Yau hypersurface $M_{k}^{k}$:
\begin{equation}
L_{m}^{k,k}(e^{t}):=1+\sum_{d=1}^{\infty}L_{m}^{k,k,d}e^{dt}\;\;(m=2,\cdots 
,k-3).
\end{equation}
\subsection{Review of Results for Fano and Calabi-Yau Hypersurfaces
and Virtual Structure Constants}
Let us summarize the results of \cite{cj}, \cite{6}. In \cite{cj}, 
we showed that the structure constants $L_{m}^{N,k,d}$ 
of $QH_{e}^{*}(M_{N}^{k})$ for $(N-k\geq 2)$  can be obtained by 
applying the recursive formulas 
which describe
$L_{m}^{N,k,d}$ in terms of $L_{m'}^{N+1,k,d'}\;\;(d'\leq d)$,
with the initial 
conditions of $L_{m}^{N,k,1}$ given by (\ref{one})
and $L_{m}^{N,k,d}=0\;\;(d\geq 2)$ in the $N\geq 2k$ region.
Let us introduce the construction of the recursive formulas given in \cite{6}.
First, we introduce the polynomial $Poly_{d}$ in 
$x,y,z_{1},z_{2},\cdots,z_{d-1}$ defined by the formula: 
\begin{eqnarray}
&&Poly_{d}(x,y,z_{1},z_{2},\cdots,z_{d-1})\no\\
&&=\frac{1}{(2\pi\sqrt{-1})^{d-1}}
\int_{C_{1}}\frac{dt_{1}}{t_{1}}\cdots
\int_{C_{d-1}}\frac{dt_{d-1}}{t_{d-1}}\prod_{j=1}^{d-1}
\biggl(\frac{(d-j)x+jy}{d}+\sum_{i=1}^{j}\frac{d-j}{d-i}t_{i}+
\sum_{i=j+1}^{d-1}\frac{j}{i}t_{i}+\no\\
&&z_{j}(\frac{(d-j)x+jy}{d}+\sum_{i=1}^{j}\frac{d-j}{d-i}t_{i}+
\sum_{i=j+1}^{d-1}\frac{j}{i}t_{i})/
(\frac{(d-j)x+jy}{d}+\sum_{i=1}^{j}\frac{d-j}{d-i}t_{i}+
\sum_{i=j+1}^{d-1}\frac{j}{i}t_{i}-z_{j})\biggr).\no\\
\label{trial2} 
\end{eqnarray}
In (\ref{trial2}), we have to choose the path $C_{i}$ carefully to
obtain the correct answer. See \cite{6} for details.
Consider the monomial 
$x^{d_{i_{0}}}z_{i_{1}}^{d_{i_{1}}}\cdots z_{i_{m}}^{d_{i_{m}}}
y^{d_{i_{m+1}}}
\;(\sum_{j=0}^{m+1}d_{i_{j}}=d-1)$, that appear in $Poly_{d}$,
associated with the following ordered partition of a positive integer $d$ 
\cite{bert}:
\begin{equation}
0=i_{0}<i_{1}<i_{2}<\cdots<i_{m}<i_{m+1}=d.
\end{equation}
Next, we prepare some elements in (a free 
abelian group) ${\bf Z}^{m+1}$, which are 
determined for each monomial 
$x^{d_{i_{0}}}z_{i_{1}}^{d_{i_{1}}}\cdots z_{i_{m}}^{d_{i_{m}}}
y^{d_{i_{m+1}}}$ ,
as follows:
\begin{eqnarray}
\alpha&:=&(m+1-d,m+1-d,\cdots,m+1-d),\no\\
\beta&:=&(0,i_{1}-1,i_{2}-2,\cdots,i_{m}-m),\no\\
\gamma&:=&(0,i_{1}(N-k),i_{2}(N-k),\cdots,i_{m}(N-k)),\no\\
\epsilon_{1}&:=&(1,0,0,0,\cdots,0),\no\\
\epsilon_{2}&:=&(1,1,0,0,\cdots,0),\no\\
\epsilon_{3}&:=&(1,1,1,0,\cdots,0),\no\\
&&\cdots\no\\
\epsilon_{m+1}&:=&(1,1,1,1,\cdots,1).
\end{eqnarray}
Now we define $\delta=(\delta_{1},\cdots,\delta_{m+1})\in{\bf Z}^{m+1}$ 
by the formula: 
\begin{equation}
\delta:=\alpha+\beta+\gamma+\sum_{j=1}^{m}(d_{i_{j}}-1)\epsilon_{j}+
d_{i_{m+1}}\epsilon_{m+1}.
\label{delta}
\end{equation}
Then the recursive formulas are given as follows:
\begin{equation}
L_{n}^{N,k,d}=\phi(Poly_{d}),
\end{equation}
where $\phi$ is a ${\bf Q}$-linear map from the ${\bf Q}$-vector
space of the homogeneous polynomials of degree $d-1$ in $x,y,z_{1},\cdots
,z_{d-1}$ to the ${\bf Q}$-vector
space of the weighted homogeneous polynomials of degree $d$ in 
$L_{m}^{N+1,k,d'}$. And it is  defined on the basis by: 
\begin{equation}
\phi(x^{d_{0}}y^{d_{d}}z_{i_{1}}^{d_{i_{1}}}\cdots z_{i_{m}}^{d_{i_{m}}})
=\prod_{j=1}^{m+1}L_{n+\delta_{j}}^{N+1,k,i_{j}-i_{j-1}}.
\end{equation}
In the $d\leq 5$ cases, we examined that 
these recursive formulas naturally lead us to the relation:
\begin{equation}
({\cal O}_{e})^{N-1}-k^{k}({\cal O}_{e})^{k-1}q=0 ,
\label{jincolli}
\end{equation}
of $QH_{e}^{*}(M_{N}^{k})\:\:(N-k\geq 2)$            
by descending induction using Beauville's result \cite{beauville},
\cite{jin}, \cite{giv}. 
In the $N-k=1$ case, the recursive formulas receive modification 
only in the $d=1$ part:
\begin{equation}
L_{m}^{k+1,k,1}=L_{m}^{k+2,k,1}-L_{0}^{k+2,k,1}=L_{m}^{k+2,k,1}-k!.
\label{givgiv}
\end{equation}
This leads us to the following relation of $QH_{e}^{*}(M_{k+1}^{k})$:
\begin{equation}
({\cal O}_{e}+k!q)^{N-1}-k^{k}({\cal O}_{e}+k!q)^{k-1}q=0.
\label{ch1}
\end{equation} 
 
The structural constant $L_{m}^{k,k,d}$ for a Calabi-Yau hypersurface 
does not obey the recursive formulas.    
Instead, we introduce here the virtual structure constants $\tilde{L}_{m}^{N,k,d}$
as follows.
\begin{defi}
Let $\tilde{L}_{m}^{N,k,d}$ be the rational number obtained by applying
the recursion relations of Fano hypersurfaces 
 for arbitrary $N$ and $k$ with the initial condition
 $L_{n}^{N,k,1}\;\;(N\geq 2k)$ and 
$L_{n}^{N,k,d}=0\;\;(d\geq 2,\;\;N\geq 2k)$. 
\end{defi}
\begin{Rem}
In the $N-k\geq2$ region, $\tilde{L}_{m}^{N,k,d}=L_{m}^{N,k,d}$.
\end{Rem}
\begin{defi}
We call $\tilde{L}_{n}^{N,k,d}$ the virtual structural constant of 
weighted degree $d$.
\end{defi}
We define the generating function of the virtual structural constants of 
the Calabi-Yau hypersurface $M_{k}^{k}$
as follows:
\begin{eqnarray}
\tilde{L}^{k,k}_{n}(e^{x})&:=&1+
\sum_{d=1}^{\infty}\tilde{L}^{k,k,d}_{n}e^{dx},\nonumber\\
&&(n=0,1,\cdots,k-1).
\label{str2}
\end{eqnarray}
In \cite{cj}, we observed that $\tilde{L}^{k,k}_{n}(e^{x})$ gives us the 
information of the B-model of the mirror manifold of $M_{k}^{k}$. 
More explicitly, we conjectured 
\begin{eqnarray}
&&\tilde{L}_{0}^{k,k}(e^{x})=
\sum_{d=0}^{\infty}\frac{(kd)!}{(d!)^{k}}e^{dx},\no\\
&&\tilde{L}_{1}^{k,k}(e^{x})=
\frac{dt(x)}{dx}:=\frac{d}{dx}
\biggl(x+(\sum_{d=1}^{\infty}\frac{(kd)!}{(d!)^{k}}\cdot 
(\sum_{i=1}^{d}\sum_{m=1}^{k-1}\frac{m}{i(ki-m)})e^{dx})/
(\sum_{d=0}^{\infty}\frac{(kd)!}{(d!)^{k}}e^{dx})\biggr)
\label{po}
\end{eqnarray}
where the r.h.s. of (\ref{po}) is derived from the solutions of 
the ODE for the period integral of the mirror manifold of $M_{k}^{k}$,
\begin{equation}
((\frac{d}{dx})^{k-1}-ke^{x}(k\frac{d}{dx}+1) 
(k\frac{d}{dx}+2)\cdots(k\frac{d}{dx}+k-1))w(x)=0,
\label{ode1}
\end{equation}
that was used in the computation based on the mirror symmetry, 
\cite{gmp}.
Of course, we can extend the conjecture (\ref{po}) to the 
general $\tilde{L}_{n}^{k,k}(e^{x})$ if we compare the 
$\tilde{L}_{n}^{k,k}(e^{x})$ with the B-model three point functions
in \cite{gmp}.
Hence we obtain the mirror map $t=t(x)$ without 
using the mirror conjecture:
\begin{equation}
t(x)=x+\int_{-\infty}^{x}dx'({\tilde{L}^{k,k}_{1}(e^{x'})}-1)
=x+\sum_{d=1}^{\infty}\frac{\tilde{L}^{k,k,d}_{1}}{d}e^{dx}.
\end{equation}
With the conjecture given by (\ref{po}), we can construct 
the mirror transformation that transforms 
the virtual structural constants of the  Calabi-Yau hypersurface into 
the real ones as follows:
\begin{equation}
L^{k,k}_{m}(e^{t})=\frac{\tilde{L}^{k,k}_{m}(e^{x(t)})}
{\tilde{L}^{k,k}_{1}(e^{x(t)})}.\quad(m=2,\cdots,k-3) 
\label{mt}
\end{equation}
The above formula is further rewritten as follows: 
\begin{equation}
L^{k,k,d}_{n}=\sum_{m=0}^{d-1}\mbox{Res}_{z=0}(z^{-m-1}
\exp(-d\sum_{j=1}^{\infty}\frac{\tilde{L}_{1}^{k,k,j}}{j}z^{j}))
\cdot(\tilde{L}^{k,k,d-m}_{n}-\tilde{L}^{k,k,d-m}_{1}).
\label{schur}
\end{equation}
This formula motivated us to propose the conjecture described in
(\ref{gene}) or in (\ref{compact}) \cite{gene}.
\section{Derivation of the Previous Results}
In this section, we first show that $V_{d-m}^{N,k,d}(n;\sigma_{m})$
satisfy the ansatz of $G_{d-m}^{N,k,d}(n;\sigma_{m})$ proposed in
\cite{gene}. Then, we show that $V_{d-m}^{N,k,d}(n;\sigma_{m})$
coincides with $G_{d-m}^{N,k,d}(n;\sigma_{m})$ in the $d\leq 3$ cases.
\begin{prop}
(i) flat metric condition
\begin{equation}
V_{d-m}^{N,k,d}(1+(k-N)d;\sigma_{m})=V_{d-m}^{N,k,d}(N-2;\sigma_{m})=0.
\end{equation}
(ii) symmetry
\begin{equation}
V_{d-m}^{N,k,d}(n;\sigma_{m})=V_{d-m}^{N,k,d}(N-1+(k-N)d-n;\sigma_{m}).
\end{equation}
(iii) 
\begin{equation}
V_{d}^{N,k,d}(n;(0))=\tilde{L}_{n}^{N,k,d}-\tilde{L}_{1+d(k-N)}^{N,k,d}.
\label{ass1}
\end{equation}
(iv)
\begin{equation}
V_{d-m}^{N,k,d}(2+(k-N)(d+f);\sigma_{m})=
V_{d-m}^{N,k,d+f}(2+(k-N)(d+f);\sigma_{m}\cup(f)).
\label{ass2}
\end{equation}
\end{prop}
{\it proof)}

(i), (ii) and (iii) are obvious by definition. 
Here, we give a proof of (iv). By the definition of 
$V_{d-m}^{N,k,d+f}(n;\sigma_{m}\cup(f))$ and by K\"ahler equation, we have
\begin{eqnarray}
&&V_{d-m}^{N,k,d+f}(2+(k-N)(d+f);\sigma_{m}\cup(f))\no\\
&&=\frac{1}{k}\frac{1}{(d-m)^{l(\sigma_{m})}}
v({\cal O}_{e^{N-2-2-(k-N)(d+f)}}{\cal O}_{e}{\cal O}_{e^{1+(k-N)f}} 
\prod_{j=1}^{l(\sigma_{m})}{\cal O}_{e^{1+(k-N)d_{j}}})_{d-m}\no\\
&&=\frac{1}{k}\frac{1}{(d-m)^{l(\sigma_{m})-1}}
v({\cal O}_{e^{N-2-2-(k-N)(d+f)}}{\cal O}_{e^{2+(k-N)(d+f)-1-(k-N)d}} 
\prod_{j=1}^{l(\sigma_{m})}{\cal O}_{e^{1+(k-N)d_{j}}})_{d-m}.\no\\
\label{reduce}
\end{eqnarray}
The last line of (\ref{reduce}) is nothing but 
$V_{d-m}^{N,k,d}(2+(k-N)(d+f);\sigma_{m})$.
 Q.E.D.

Then we introduce the definition:
\begin{defi}
\begin{equation}
\tilde{V}_{d-m}^{N,k,d+f}(n;\sigma_{m}\cup(f)):=
\sum_{j=0}^{(k-N)f}V_{d-m}^{N,k,d}(n-j;\sigma_{m})-
\sum_{j=0}^{(k-N)f}V_{d-m}^{N,k,d}(1+(k-N)(d+f)-j;\sigma_{m}).
\label{gnew}
\end{equation}
We denote by $\pi_{f}$ the map which maps the function 
$g(n)$ on ${\bf Z}$ to\\ $\sum_{j=0}^{(k-N)f}g(n-j)-
\sum_{j=0}^{(k-N)f}g(1+(k-N)(d+f)-j)$. Then we have
\begin{equation}
\pi_{f}(V_{d-m}^{N,k,d}(n;\sigma_{m}))=
\tilde{V}_{d-m}^{N,k,d+f}(n;\sigma_{m}\cup(f)).
\end{equation}
\end{defi}
With this definition, we are led to the following proposition.
\begin{prop}
$\tilde{V}_{d-m}^{N,k,d+f}(n;\sigma_{m}\cup(f))$
satisfies the condition (iv) of the above proposition:
\begin{equation}
\tilde{V}_{d-m}^{N,k,d+f}(2+(k-N)(d+f);\sigma_{m}\cup(f))
=V_{d-m}^{N,k,d}(2+(k-N)(d+f);\sigma_{m}),
\label{sol}
\end{equation}
and the condition (i) and (ii).
\end{prop}
{\it proof)} 
The fact that $\tilde{V}_{d-m}^{N,k,d+f}(n;\sigma_{m}\cup(f))$
satisfies the condition (i) is rather obvious,
and we first prove the condition (ii).
It suffices to consider the part
$\sum_{j=0}^{(k-N)f}V_{d-m}^{N,k,d}(n-j;\sigma_{m})$:
\begin{eqnarray}
&&\sum_{j=0}^{(k-N)f}V_{d-m}^{N,k,d}(N-1+(k-N)(d+f)-n-j;\sigma_{m})\no\\
&=&\sum_{j=0}^{(k-N)f}V_{d-m}^{N,k,d}(N-1+(k-N)d-n+j;\sigma_{m})\no\\
&=&\sum_{j=0}^{(k-N)f}V_{d-m}^{N,k,d}(n-j;\sigma_{m}).
\label{(ii)}
\end{eqnarray}  
Next, we turn to the formula 
$(\ref{sol})$. By definition and by the condition
$V_{d-m}^{N,k,d}(1+(k-N)d;\sigma_{m})=0$, we obtain 
\begin{eqnarray}
&&\tilde{V}_{d-m}^{N,k,d+f}(2+(k-N)(d+f);\sigma_{m}\cup(f))\no\\
&=&\sum_{j=0}^{(k-N)f}V_{d-m}^{N,k,d}(2+(k-N)(d+f)-j;\sigma_{m})
-\sum_{j=0}^{(k-N)f}V_{d-m}^{N,k,d}(1+(k-N)(d+f)-j;\sigma_{m})\no\\
&=&V_{d-m}^{N,k,d}(2+(k-N)(d+f);\sigma_{m}).
\label{prf}
\end{eqnarray}
 Q.E.D.
\begin{defi}
\begin{equation}
hi_{d-m}^{N,k,d+f}(n;\sigma_{m}\cup(f)):=
V_{d-m}^{N,k,d+f}(n;\sigma_{m}\cup(f))-
\tilde{V}_{d-m}^{N,k,d+f}(n;\sigma_{m}\cup(f))
\label{difference}
\end{equation}
\end{defi}
\begin{prop}
$hi_{d-m}^{N,k,d+f}(n;\sigma_{m}\cup(f))$ satisfies 
the condition:
\begin{equation}
hi_{d-m}^{N,k,d+f}(2+(k-N)(d+f);\sigma_{m}\cup(f))
=hi_{d-m}^{N,k,d+f}(1+(k-N)(d+f);\sigma_{m}\cup(f))=0,
\label{hidden}
\end{equation}
and the condition (i), (ii).
\end{prop}
{\it proof)} Immediate. Q.E.D.
\\
With the above discussions, we can construct the following decomposition 
of $V_{d-m}^{N,k,d}(n;\sigma_{m})\quad (\sigma_{m}:m=d_{1}+d_{2}+\cdots+
d_{l(\sigma_{m})},\;d_{1}\geq\cdots\geq d_{l(\sigma_{m})}\geq1)$:
\begin{eqnarray}
V_{d-m}^{N,k,d}(n;\sigma_{m})&=&\pi_{d_{l(\sigma_{m})}}
\circ\cdots\circ\pi_{d_{2}}\circ\pi_{d_{1}}
(\tilde{L}^{N,k,d}_{n}-\tilde{L}^{N,k,d}_{1+(k-N)d})\no\\
&&+\sum_{j=2}^{l(\sigma_{m})}\pi_{d_{l(\sigma_{m})}}
\circ\cdots\circ\pi_{d_{j}}(hi^{(j)})+hi^{(l(\sigma_{m}))}.
\label{dec}
\end{eqnarray}
\begin{Rem}
The above decomposition depends on the order of $d_{i}$'s. Hence it is not
 unique.
\end{Rem}
In (\ref{dec}), we omit the subscripts of 
$hi_{d-m}^{N,k,d+f}(n;\sigma_{m}\cup(f))$.   
We can easily see that $hi_{d-m}^{N,k,d+f}(n;\sigma_{m}\cup(f))$
consists of monomials of degree $(d-m)$ of 
$\tilde{L}_{j}^{N,k,m'}\;\;(m'<d-m)$ only, because the linear dependence 
on $\tilde{L}_{j}^{N,k,d-m}$ cannot satisfy the condition (\ref{hidden}).
Thus we are led to the
proposition by picking up the top term of the decomposition in (\ref{dec}):
\begin{prop}
The linear part of $V_{d-m}^{N,k,d}(n;\sigma_{m})$ is given by the formula:
\begin{eqnarray}
&&\sum_{j_{1}=0}^{d_{1}(k-N)}
\sum_{j_{2}=0}^{d_{2}(k-N)}\cdots
\sum_{j_{l(\sigma_{m})}=0}^{d_{l(\sigma_{m})}(k-N)}
(\tilde{L}^{N,k,d-m}_{n-\sum_{i=1}^{l(\sigma_{m})}j_{i}}-\tilde{L}^{N,k,d-m}_{1+(k-N)d-\sum_{i=1}^{l(\sigma_{m})}j_{i}}).\no\\
\label{genetop1}
\end{eqnarray}
Especially in the $d-m=1$ case, we obtain the following equality:
\begin{eqnarray}
V_{1}^{N,k,d}(n;\sigma_{d-1})&=&\sum_{j_{1}=0}^{d_{1}(k-N)}
\sum_{j_{2}=0}^{d_{2}(k-N)}\cdots
\sum_{j_{l(\sigma_{d-1})}=0}^{d_{l(\sigma_{d-1})}(k-N)}
(\tilde{L}^{N,k,1}_{n-\sum_{i=1}^{l(\sigma_{d-1})}j_{i}}
-\tilde{L}^{N,k,1}_{1+(k-N)d-\sum_{i=1}^{l(\sigma_{d-1})}j_{i}}).\no\\
\label{genetop2}
\end{eqnarray}
\end{prop}
\begin{Rem}
By introducing the polynomial in $x$:
\begin{eqnarray}
\prod_{j=1}^{l(\sigma_{m})}\frac{(1-x^{d_{j}(k-N)+1})}{(1-x)}=
\sum_{j=0}^{(k-N)m}A_{j}^{N,k}(\sigma_{m})x^{j},
\end{eqnarray}
we can write the linear part (\ref{genetop1}) in a more compact form,
\begin{eqnarray}
&&\sum_{j_{1}=0}^{d_{1}(k-N)}
\sum_{j_{2}=0}^{d_{2}(k-N)}\cdots
\sum_{j_{l(\sigma_{m})}=0}^{d_{l(\sigma_{m})}(k-N)}
(\tilde{L}^{N,k,d-m}_{n-\sum_{i=1}^{l(\sigma_{m})}j_{i}}-
\tilde{L}^{N,k,d-m}_{1+(k-N)d-\sum_{i=1}^{l(\sigma_{m})}j_{i}})\no\\
&&=\sum_{j=0}^{(k-N)m}A_{j}^{N,k}(\sigma_{m})(\tilde{L}_{n-j}^{N,k,d-m}
-\tilde{L}_{1+(k-N)d-j}^{N,k,d-m}).
\end{eqnarray}
\end{Rem}
These are the results directly applying the constraints obtained in 
Proposition 1. 
On the other hand, we can explicitly
express $V_{d-m}^{N,k,d}(n;\sigma_{m})$ in terms of 
$\tilde{L}_{n}^{N,k,d}$, because the virtual Gromov-Witten invariants 
satisfy the K\"ahler equation and the associativity equation. 
As examples, we compute $V_{d-m}^{N,k,d}(n;\sigma_{m})\;\;(d\leq 3)$.    
\begin{prop}
$V_{d-m}^{N,k,d}(n;\sigma_{m})\;\;(d\leq 3)$ can be written in terms of
 $\tilde{L}^{N,k,d}_{n}\;\;(d\leq 3)$ as follows.
\begin{eqnarray}
V_{1}^{N,k,1}(n;(0))&=&\tilde{L}_{n}^{N,k,1}-\tilde{L}_{1+(k-N)}^{N,k,1},\no\\
V_{2}^{N,k,2}(n;(0))&=&\tilde{L}_{n}^{N,k,2}-\tilde{L}_{1+2(k-N)}^{N,k,2},\no\\
V_{1}^{N,k,2}(n;(1))&=&\sum_{j=0}^{k-N}
(\tilde{L}_{n-j}^{N,k,1}-\tilde{L}_{1+2(k-N)-j}^{N,k,1}),
\no\\
V_{3}^{N,k,3}(n;(0))&=&\tilde{L}_{n}^{N,k,3}-\tilde{L}_{1+3(k-N)}^{N,k,3},
\no\\
V_{2}^{N,k,3}(n;(1))&=&\sum_{j=0}^{k-N}
(\tilde{L}_{n-j}^{N,k,2}-\tilde{L}_{1+3(k-N)-j}^{N,k,2})+hi_{2}^{N,k,3}
(n;(1))
\no\\
V_{1}^{N,k,3}(n;(2))&=&
\sum_{j=0}^{2(k-N)}(\tilde{L}_{n-j}^{N,k,1}
-\tilde{L}_{1+3(k-N)-j}^{N,k,1}),\no\\
V_{1}^{N,k,3}(n;(1)+(1))&=&
\sum_{j=0}^{2(k-N)}A_{j}^{N,k}((1)+(1))(\tilde{L}_{n-j}^{N,k,1}
-\tilde{L}_{1+3(k-N)-j}^{N,k,1}).
\label{cubic}
\end{eqnarray}
where $hi_{2}^{N,k,3}(n;(1))$ is
given by the formula:
\begin{eqnarray}
hi_{2}^{N,k,3}(n;(1))&=&
\sum_{j=0}^{(k-N)-1}(\sum_{m=0}^{j}\tilde{L}_{n-m}^{N,k,1}
\tilde{L}_{n-2(k-N)+j-m}^{N,k,1}-\tilde{L}_{(k-N)+2+j}^{N,k,1}
\sum_{m=0}^{2(k-N)}\tilde{L}_{n-m}^{N,k,1}\no\\
&&+\tilde{L}_{1+(k-N)}^{N,k,1}
\sum_{m=j+1}^{2(k-N)-j-1}\tilde{L}_{n-m}^{N,k,1})\no\\
&&-\sum_{j=0}^{(k-N)-1}(\sum_{m=0}^{j}\tilde{L}_{1+3(k-N)-m}^{N,k,1}
\tilde{L}_{1+(k-N)+j-m}^{N,k,1}-\tilde{L}_{(k-N)+2+j}^{N,k,1}
\sum_{m=0}^{2(k-N)}\tilde{L}_{1+3(k-N)-m}^{N,k,1}\no\\
&&+\tilde{L}_{1+(k-N)}^{N,k,1}
\sum_{m=j+1}^{2(k-N)-j-1}\tilde{L}_{1+3(k-N)-m}^{N,k,1}).
\end{eqnarray}
\end{prop}
{\it proof)} 

We give a proof of the formula for $V_{2}^{N,k,3}(n;(1))$, because 
the other formulas follow
obviously from the preceding discussions. First, we introduce the 
following virtual G-W invariant:
\begin{eqnarray}
v({\cal O}_{e^{N-2-n}}{\cal O}_{e^{n-1-2(k-N)-m}}{\cal O}_{e^{1+m}})_{2}.
\end{eqnarray}
Using the associativity equation, we obtain the equality, 
\begin{eqnarray}
&&\frac{1}{k}\cdot 
v({\cal O}_{e^{N-2-n}}{\cal O}_{e^{n-1-2(k-N)-m}}{\cal O}_{e^{1+m}})_{2}
-\frac{1}{k}\cdot 
v({\cal O}_{e^{N-2-n}}{\cal O}_{e^{n-2(k-N)-m}}{\cal O}_{e^{m}})_{2}\no\\
&&=\tilde{L}_{n-m}^{N,k,2}-\tilde{L}_{1+m+2(k-N)}^{N,k,2}\no\\
&&+\sum_{j=0}^{m-1}(\tilde{L}_{n-j}^{N,k,1}-\tilde{L}_{1+(k-N)+j}^{N,k,1})
\cdot(\tilde{L}_{n-m-(k-N)}^{N,k,1}-\tilde{L}_{1+(k-N)}^{N,k,1})\no\\
&&-\sum_{j=0}^{m+(k-N)}(\tilde{L}_{n-j}^{N,k,1}-\tilde{L}_{1+(k-N)+j}^{N,k,1})
\cdot(\tilde{L}_{1+m+(k-N)}^{N,k,1}-\tilde{L}_{1+(k-N)}^{N,k,1}). 
\label{conic}
\end{eqnarray}
Adding up (\ref{conic}) with $m=0,1,\cdots,(k-N)$ and using some
algebras, we can reach the desired formula in the proposition. Q.E.D.

From this proposition and the conjectural form of the generalized mirror 
transformation for the $d\leq 3$ rational G-W invariants in \cite{gene}, 
we can see the equality 
$G_{d-m}^{N,k,d}(n;\sigma_{m})=V_{d-m}^{N,k,d}(n;\sigma_{m})\;\;(d\leq 3)$.
Thus, we have derived the generalized mirror transformation up to $d\leq 
3$ cases in \cite{gene} under the assumption of Conjecture 2.
\section{Explicit Determination in the $k-N=1,\;d=4,5$ cases}
In this section, we restrict $k-N$ to $1$ to avoid the complication of
the formulas. 
In this setting, we compute $V_{d-m}^{k-1,k,d}(n;\sigma_{m})$ in the
$d=4,5$ cases and discuss the modification of 
$V_{d-m}^{k-1,k,d}(n;\sigma_{m})$ into $G_{d-m}^{k-1,k,d}(n;\sigma_{m})$.
With the aid of some numerical data, we fix the modification and 
derive the generalized mirror transformation in these
cases. Generalization to the general $k-N$ cases is rather straightforward. 

First, we repeatedly use the associativity equation and obtain the 
following formula that represent $V_{4-m}^{k-1,k,4}(n;\sigma_{m})$ in
terms of $\tilde{L}^{k-1,k,d}_{n}$. 
\begin{prop}
$V_{4-m}^{k-1,k,4}(n;\sigma_{m})$'s are inductively determined as follows:
\begin{eqnarray}
V_{4}^{k-1,k,4}(n;(0))&=&\tilde{L}^{k-1,k,4}_{n}-\tilde{L}^{k-1,k,4}_{5},\no\\
V_{3}^{k-1,k,4}(n;(1))&=&\tilde{L}^{k-1,k,3}_{n}+\tilde{L}^{k-1,k,3}_{n-1}-
\tilde{L}^{k-1,k,3}_{5}-\tilde{L}^{k-1,k,3}_{4}\no\\
&&+(\tilde{L}_{n}^{k-1,k,2}-\tilde{L}_{3}^{k-1,k,2})\cdot
(\tilde{L}_{n-3}^{k-1,k,1}-\tilde{L}_{2}^{k-1,k,1})\no\\
&&+(\tilde{L}_{n}^{k-1,k,1}-\tilde{L}_{2}^{k-1,k,1})\cdot
(\tilde{L}_{n-2}^{k-1,k,2}-\tilde{L}_{3}^{k-1,k,2})\no\\
&&-(\tilde{L}_{3}^{k-1,k,1}-\tilde{L}_{2}^{k-1,k,1})\cdot 
V_{2}^{k-1,k,4}(n;(2))\no\\
&&-V_{1}^{k-1,k,4}(n;(3))\cdot 
(\tilde{L}_{4}^{k-1,k,2}-\tilde{L}_{3}^{k-1,k,2}),\no\\
V_{2}^{k-1,k,4}(n;(2))&=&V_{2}^{k-1,k,3}(n;(1))+
\tilde{L}^{k-1,k,2}_{n-2}-\tilde{L}^{k-1,k,2}_{5}\no\\
&&+V_{1}^{k-1,k,3}(n;(2))
\cdot(\tilde{L}_{n-3}^{k-1,k,1}-\tilde{L}_{2}^{k-1,k,1})\no\\
&&-V_{1}^{k-1,k,4}(n;(3))\cdot 
(\tilde{L}_{4}^{k-1,k,1}-\tilde{L}_{2}^{k-1,k,1}),\no\\
V_{1}^{k-1,k,4}(n;(3))&=&\tilde{L}^{k-1,k,1}_{n}+\tilde{L}^{k-1,k,1}_{n-1}+
\tilde{L}^{k-1,k,1}_{n-2}+\tilde{L}^{k-1,k,1}_{n-3}\no\\
&&-(\tilde{L}^{k-1,k,1}_{5}+
\tilde{L}^{k-1,k,2}_{4}+\tilde{L}^{k-1,k,2}_{3}+\tilde{L}^{k-1,k,1}_{2}),\no\\
V_{1}^{k-1,k,4}(n;(1)+(2))&=&\tilde{L}^{k-1,k,1}_{n}+2
\tilde{L}^{k-1,k,1}_{n-1}+
2\tilde{L}^{k-1,k,1}_{n-2}+\tilde{L}^{k-1,k,1}_{n-3}\no\\
&&-(\tilde{L}^{k-1,k,1}_{5}+
2\tilde{L}^{k-1,k,2}_{4}+2\tilde{L}^{k-1,k,2}_{3}+\tilde{L}^{k-1,k,1}_{2}),
\no\\
V_{1}^{k-1,k,4}(n;(1)+(1)+(1))&=&
\tilde{L}^{k-1,k,1}_{n}+3\tilde{L}^{k-1,k,1}_{n-1}+
3\tilde{L}^{k-1,k,1}_{n-2}+\tilde{L}^{k-1,k,1}_{n-3}\no\\
&&-(\tilde{L}^{k-1,k,1}_{5}+
3\tilde{L}^{k-1,k,2}_{4}+3\tilde{L}^{k-1,k,2}_{3}+
\tilde{L}^{k-1,k,1}_{2}),\no\\
V_{2}^{k-1,k,4}(n;(1)+(1))&=&V_{2}^{k-1,k,3}(n;(1))+V_{2}^{k-1,k,3}(n-1;(1))
-(\tilde{L}_{5}^{k-1,k,2}-\tilde{L}_{3}^{k-1,k,2})\no\\
&&+\frac{1}{2}\biggl(V_{1}^{k-1,k,1}(n;(1))\cdot
(\tilde{L}_{n-3}^{k-1,k,1}-\tilde{L}_{2}^{k-1,k,1})\no\\
&&+(\tilde{L}_{n}^{k-1,k,1}-\tilde{L}_{2}^{k-1,k,1})\cdot
V_{1}^{k-1,k,1}(n-2;(1))\no\\
&&-V_{1}^{k-1,k,4}(n;(1)+(2))\cdot
(\tilde{L}_{3}^{k-1,k,1}-\tilde{L}_{2}^{k-1,k,1})\no\\
&&-V_{1}^{k-1,k,4}(n;(3))\cdot
(\tilde{L}_{4}^{k-1,k,1}-\tilde{L}_{2}^{k-1,k,1})\biggr).  
\label{4}
\end{eqnarray}
\end{prop}
\begin{conj}
In the  $d=4$ case, $V_{4-m}^{k-1,k,4}(n;\sigma_{m})$ equals  
$G_{4-m}^{k-1,k,4}(n;\sigma_{m})$ except for
$V_{2}^{k-1,k,4}(n,(1)+(1))$. $G_{2}^{k-1,k,4}(n,(1)+(1))$ is given by
the formula:
\begin{eqnarray}
G_{2}^{k-1,k,4}(n;(1)+(1))&=&V_{2}^{k-1,k,3}(n;(1))+V_{2}^{k-1,k,3}(n-1;(1))
-(\tilde{L}_{5}^{k-1,k,2}-\tilde{L}_{3}^{k-1,k,2})\no\\
&&+\frac{3}{4}\biggl(V_{1}^{k-1,k,1}(n;(1))\cdot
(\tilde{L}_{n-3}^{k-1,k,1}-\tilde{L}_{2}^{k-1,k,1})\no\\
&&+(\tilde{L}_{n}^{k-1,k,1}-\tilde{L}_{2}^{k-1,k,1})\cdot
V_{1}^{k-1,k,1}(n-2;(1))\no\\
&&-V_{1}^{k-1,k,4}(n;(1)+(2))\cdot
(\tilde{L}_{3}^{k-1,k,1}-\tilde{L}_{2}^{k-1,k,1})\no\\
&&-V_{1}^{k-1,k,4}(n;(3))\cdot
(\tilde{L}_{4}^{k-1,k,1}-\tilde{L}_{2}^{k-1,k,1})\biggr).
\label{modify}
\end{eqnarray}
\end{conj}
The modification of the factor $\frac{1}{2}$ in (\ref{4}) into 
the factor $\frac{3}{4}$ in (\ref{modify}) is determined by one
numerical data:
\begin{equation}
L_{7}^{11,12,4}=1324882975682876246483412831870565329165165953902032,
\label{data1}
\end{equation}
and the generalized mirror transformation obtained by 
Conjecture 3 and (\ref{gene}) correctly predicts $L_{n}^{k-1,k,4}$ 
up to $k\leq 18$.

In the $d=5$ case, we can state the following proposition under the
assumption of Conjecture 2.
\begin{prop}
\begin{eqnarray}
&&G_{5}^{k-1,k,5}(n;(0))=V_{5}^{k-1,k,5}(n;(0)),\no\\
&&G_{4}^{k-1,k,5}(n;(1))=V_{4}^{k-1,k,5}(n;(1)),\no\\
&&G_{3}^{k-1,k,5}(n;(2))=V_{3}^{k-1,k,5}(n;(2)),\no\\
&&G_{2}^{k-1,k,5}(n;(3))=V_{2}^{k-1,k,5}(n;(3)),\no\\
&&G_{1}^{k-1,k,5}(n;(4))=V_{1}^{k-1,k,5}(n;(4)),\no\\
&&G_{1}^{k-1,k,5}(n;(3)+(1))=V_{1}^{k-1,k,5}(n;(3)+(1)),\no\\
&&G_{1}^{k-1,k,5}(n;(2)+(2))=V_{1}^{k-1,k,5}(n;(2)+(2)),\no\\
&&G_{1}^{k-1,k,5}(n;(2)+(1)+(1))=V_{1}^{k-1,k,5}(n;(2)+(1)+(1)),\no\\
&&G_{1}^{k-1,k,5}(n;(1)+(1)+(1)+(1))=V_{1}^{k-1,k,5}(n;(1)+(1)+(1)+(1)).
\end{eqnarray}
\end{prop} 
On the other hand, we find that there exist some non-trivial modifications 
to obtain $G_{3}^{k-1,k,5}(n;(1)+(1))$, $G_{2}^{k-1,k,5}(n;(1)+(2))$ and 
$G_{2}^{k-1,k,5}(n;(1)+(1)+(1))$ from the corresponding virtual Gromov-
Witten Invariants. 
\begin{conj}
\begin{eqnarray}
G_{2}^{k-1,k,5}(n;(1)+(2))&=&V_{2}^{k-1,k,4}(n;(2))+V_{2}^{k-1,k,4}(n-1;(2))
-\tilde{L}^{k-1,k,2}_{6}+
\tilde{L}^{k-1,k,2}_{3}\no\\
&&+\frac{8}{5}hi_{1}(n)+hi_{2}(n)+\frac{4}{5}hi_{3}(n)-\frac{3}{5}hi_{4}(n),
\end{eqnarray}
\begin{eqnarray}
G_{2}^{k-1,k,5}(n;(1)+(1)+(1))&=& 
V_{2}^{k-1,k,3}(n;(1))+V_{2}^{k-1,k,3}(n-1;(1))
-(\tilde{L}_{5}^{k-1,k,2}-\tilde{L}_{3}^{k-1,k,2})\no\\
&&+\frac{4}{5}\biggl(V_{1}^{k-1,k,1}(n;(1))\cdot
(\tilde{L}_{n-3}^{k-1,k,1}-\tilde{L}_{2}^{k-1,k,1})\no\\
&&+(\tilde{L}_{n}^{k-1,k,1}-\tilde{L}_{2}^{k-1,k,1})\cdot
V_{1}^{k-1,k,1}(n-2;(1))\no\\
&&-V_{1}^{k-1,k,4}(n;(1)+(2))\cdot
(\tilde{L}_{3}^{k-1,k,1}-\tilde{L}_{2}^{k-1,k,1})\no\\
&&-V_{1}^{k-1,k,4}(n;(3))\cdot
(\tilde{L}_{4}^{k-1,k,1}-\tilde{L}_{2}^{k-1,k,1})\biggr)\no\\
&&+V_{2}^{k-1,k,3}(n-1;(1))+V_{2}^{k-1,k,3}(n-2;(1))
-(\tilde{L}_{5}^{k-1,k,2}-\tilde{L}_{3}^{k-1,k,2})\no\\
&&+\frac{4}{5}\biggl(V_{1}^{k-1,k,1}(n-1;(1))\cdot
(\tilde{L}_{n-4}^{k-1,k,1}-\tilde{L}_{2}^{k-1,k,1})\no\\
&&+(\tilde{L}_{n-1}^{k-1,k,1}-\tilde{L}_{2}^{k-1,k,1})\cdot
V_{1}^{k-1,k,1}(n-3;(1))\no\\
&&-V_{1}^{k-1,k,4}(n-1;(1)+(2))\cdot
(\tilde{L}_{3}^{k-1,k,1}-\tilde{L}_{2}^{k-1,k,1})\no\\
&&-V_{1}^{k-1,k,4}(n-1;(3))\cdot
(\tilde{L}_{4}^{k-1,k,1}-\tilde{L}_{2}^{k-1,k,1})\biggr)
-V_{2}^{k-1,k,3}(6;(1))\no\\
&&+\frac{46}{25}hi_{1}(n)+\frac{46}{25}hi_{2}(n)+\frac{16}{25}hi_{3}(n)-
\frac{2}{25}hi_{4}(n),
\end{eqnarray}
\begin{eqnarray}
G_{3}^{k-1,k,5}(n;(1)+(1))&=&V_{3}^{k-1,k,4}(n;(1))+V_{3}^{k-1,k,4}(n-1;(1)) 
-(\tilde{L}_{6}^{k-1,k,3}-
\tilde{L}_{4}^{k-1,k,3})\no\\
&&+\frac{4}{5}\biggl(V_{2}^{k-1,k,3}(n;(1))\cdot(\tilde{L}_{n-4}^{k-1,k,1}-
\tilde{L}_{2}^{k-1,k,1})\no\\
&&+(\tilde{L}_{n}^{k-1,k,1}-
\tilde{L}_{2}^{k-1,k,1})\cdot V_{2}^{k-1,k,3}(n-2;(1))\no\\
&&-\bigl(V_{2}^{k-1,k,4}(n;(2))+V_{2}^{k-1,k,4}(n-1;(2))\no\\
&&-(\tilde{L}^{k-1,k,2}_{6}-
\tilde{L}^{k-1,k,2}_{3})\bigr)\cdot(\tilde{L}_{3}^{k-1,k,1}-
\tilde{L}_{2}^{k-1,k,1})\no\\
&&-V_{1}^{k-1,k,5}(n;(4))\cdot(\tilde{L}_{5}^{k-1,k,2}-
\tilde{L}_{3}^{k-1,k,2})\biggr)\no\\
&&+\frac{3}{5}\biggl(V_{1}^{k-1,k,2}(n;(1))\cdot(\tilde{L}_{n-3}^{k-1,k,2}-
\tilde{L}_{3}^{k-1,k,2})\no\\
&&+(\tilde{L}_{n}^{k-1,k,2}-
\tilde{L}_{3}^{k-1,k,2})\cdot V_{1}^{k-1,k,2}(n-3;(1))\no\\
&&-V_{1}^{k-1,k,5}(n;(1)+(3))\cdot
(\tilde{L}_{4}^{k-1,k,2}-\tilde{L}_{3}^{k-1,k,2})\no\\
&&-V_{2}^{k-1,k,5}(n;(3))
\cdot(\tilde{L}_{4}^{k-1,k,1}-\tilde{L}_{2}^{k-1,k,1})\biggr)\no\\
&&-(\tilde{L}_{3}^{k-1,k,1}-\tilde{L}_{2}^{k-1,k,1})\cdot
(\frac{6}{5}hi_{1}(n)+hi_{2}(n)+\frac{3}{5}hi_{3}(n)-
\frac{1}{5}hi_{4}(n)),
\no\\     
\end{eqnarray}
where $hi_{j}(n)$ is a degree $2$ homogeneous polynomial of 
$\tilde{L}^{k-1,k,1}_{m}$ satisfying $hi_{j}(6)=hi_{j}(7)=0$ and is
given by,  
\begin{eqnarray}
hi_{1}(n)&=&\tilde{L}_{n}^{k-1,k,1}\tilde{L}_{n-4}^{k-1,k,1}-
\tilde{L}_{3}^{k-1,k,1}(\tilde{L}_{n}^{k-1,k,1}+\tilde{L}_{n-1}^{k-1,k,1}+
\tilde{L}_{n-2}^{k-1,k,1}+\tilde{L}_{n-3}^{k-1,k,1}+\tilde{L}_{n-4}^{k-1,k,1})
\no\\
&&+\tilde{L}_{2}^{k-1,k,1}(\tilde{L}_{n-1}^{k-1,k,1}+\tilde{L}_{n-2}^{k-1,k,1}+\tilde{L}_{n-3}^{k-1,k,1})\no\\
&&-(\tilde{L}_{6}^{k-1,k,1}\tilde{L}_{2}^{k-1,k,1}-
\tilde{L}_{3}^{k-1,k,1}(\tilde{L}_{6}^{k-1,k,1}+\tilde{L}_{5}^{k-1,k,1}+
\tilde{L}_{4}^{k-1,k,1}+\tilde{L}_{3}^{k-1,k,1}+\tilde{L}_{2}^{k-1,k,1})\no\\
&&+\tilde{L}_{2}^{k-1,k,1}(\tilde{L}_{5}^{k-1,k,1}+
\tilde{L}_{4}^{k-1,k,1}+\tilde{L}_{3}^{k-1,k,1})),\no\\
hi_{2}(n)&=&\tilde{L}_{n}^{k-1,k,1}\tilde{L}_{n-3}^{k-1,k,1}
+\tilde{L}_{n-1}^{k-1,k,1}\tilde{L}_{n-4}^{k-1,k,1}\no\\
&&-\tilde{L}_{4}^{k-1,k,1}(\tilde{L}_{n}^{k-1,k,1}+\tilde{L}_{n-1}^{k-1,k,1}+
\tilde{L}_{n-2}^{k-1,k,1}+\tilde{L}_{n-3}^{k-1,k,1}+\tilde{L}_{n-4}^{k-1,k,1})
+\tilde{L}_{2}^{k-1,k,1}\tilde{L}_{n-2}^{k-1,k,1}\no\\
&&-(\tilde{L}_{6}^{k-1,k,1}\tilde{L}_{3}^{k-1,k,1}
+\tilde{L}_{5}^{k-1,k,1}\tilde{L}_{2}^{k-1,k,1}\no\\
&&-\tilde{L}_{4}^{k-1,k,1}(\tilde{L}_{6}^{k-1,k,1}+\tilde{L}_{5}^{k-1,k,1}+
\tilde{L}_{4}^{k-1,k,1}+\tilde{L}_{3}^{k-1,k,1}+\tilde{L}_{2}^{k-1,k,1})
+\tilde{L}_{2}^{k-1,k,1}\tilde{L}_{4}^{k-1,k,1}),\no\\
hi_{3}(n)&=&\tilde{L}_{n}^{k-1,k,1}\tilde{L}_{n-2}^{k-1,k,1}
+\tilde{L}_{n-1}^{k-1,k,1}\tilde{L}_{n-3}^{k-1,k,1}+
\tilde{L}_{n-2}^{k-1,k,1}\tilde{L}_{n-4}^{k-1,k,1}\no\\
&&-\tilde{L}_{5}^{k-1,k,1}(\tilde{L}_{n}^{k-1,k,1}+\tilde{L}_{n-1}^{k-1,k,1}+
\tilde{L}_{n-2}^{k-1,k,1}+\tilde{L}_{n-3}^{k-1,k,1}+\tilde{L}_{n-4}^{k-1,k,1})
-\tilde{L}_{2}^{k-1,k,1}\tilde{L}_{n-2}^{k-1,k,1}\no\\
&&-(\tilde{L}_{6}^{k-1,k,1}\tilde{L}_{4}^{k-1,k,1}
+\tilde{L}_{5}^{k-1,k,1}\tilde{L}_{3}^{k-1,k,1}+
\tilde{L}_{4}^{k-1,k,1}\tilde{L}_{2}^{k-1,k,1}\no\\
&&-\tilde{L}_{5}^{k-1,k,1}(\tilde{L}_{6}^{k-1,k,1}+\tilde{L}_{5}^{k-1,k,1}+
\tilde{L}_{4}^{k-1,k,1}+\tilde{L}_{3}^{k-1,k,1}+\tilde{L}_{2}^{k-1,k,1})
-\tilde{L}_{2}^{k-1,k,1}\tilde{L}_{4}^{k-1,k,1}),\no\\
hi_{4}(n)&=&\tilde{L}_{n-1}^{k-1,k,1}\tilde{L}_{n-3}^{k-1,k,1}
-\tilde{L}_{4}^{k-1,k,1}(\tilde{L}_{n-1}^{k-1,k,1}+
\tilde{L}_{n-2}^{k-1,k,1}+\tilde{L}_{n-3}^{k-1,k,1})
+\tilde{L}_{3}^{k-1,k,1}\tilde{L}_{n-2}^{k-1,k,1}\no\\
&&-(\tilde{L}_{5}^{k-1,k,1}\tilde{L}_{3}^{k-1,k,1}
-\tilde{L}_{4}^{k-1,k,1}(\tilde{L}_{5}^{k-1,k,1}+
\tilde{L}_{4}^{k-1,k,1}+\tilde{L}_{3}^{k-1,k,1})
+\tilde{L}_{3}^{k-1,k,1}\tilde{L}_{4}^{k-1,k,1}).\no\\
\end{eqnarray}
\end{conj}
As in the $d=4$ case, we have fixed the  modification of the rational
factors by one numerical data:
\begin{eqnarray}
&&L_{8}^{12,13,5}=\no\\
&&100355724573836807695163109854598526931747042477505803923089934593470758513921/180000.\no\\
\end{eqnarray}
Check of the prediction formula obtained by (\ref{gene}), Proposition 7
and Conjecture 4 takes a lot of time due to numerical computation of
fixed point formulas. We checked that the formula correctly predicts 
$L_{8}^{13,14,5}$.
Now, we discuss the rules of modification of $V_{d-m}^{k-1,k,d}(n;\sigma_{m})$
into $G_{d-m}^{k-1,k,d}(n;\sigma_{m})$.
First, the corresponding $V_{5-m}^{k-1,k,5}(n;\sigma_{m})$ is given as follows.
 \begin{eqnarray}
V_{2}^{k-1,k,5}(n;(1)+(2))&=&V_{2}^{k-1,k,4}(n;(2))+V_{2}^{k-1,k,4}(n-1;(2))
-\tilde{L}^{k-1,k,2}_{6}+
\tilde{L}^{k-1,k,2}_{3}\no\\
&&+\frac{1}{2}(2hi_{1}(n)+hi_{2}(n)+hi_{3}(n)-hi_{4}(n)),
\end{eqnarray}
\begin{eqnarray}
V_{2}^{k-1,k,5}(n;(1)+(1)+(1))&=& 
V_{2}^{k-1,k,3}(n;(1))+V_{2}^{k-1,k,3}(n-1;(1))
-(\tilde{L}_{5}^{k-1,k,2}-\tilde{L}_{3}^{k-1,k,2})\no\\
&&+\frac{1}{2}\biggl(V_{1}^{k-1,k,1}(n;(1))\cdot
(\tilde{L}_{n-3}^{k-1,k,1}-\tilde{L}_{2}^{k-1,k,1})\no\\
&&+(\tilde{L}_{n}^{k-1,k,1}-\tilde{L}_{2}^{k-1,k,1})\cdot
V_{1}^{k-1,k,1}(n-2;(1))\no\\
&&-V_{1}^{k-1,k,4}(n;(1)+(2))\cdot
(\tilde{L}_{3}^{k-1,k,1}-\tilde{L}_{2}^{k-1,k,1})\no\\
&&-V_{1}^{k-1,k,4}(n;(3))\cdot
(\tilde{L}_{4}^{k-1,k,1}-\tilde{L}_{2}^{k-1,k,1})\biggr)\no\\
&&+V_{2}^{k-1,k,3}(n-1;(1))+V_{2}^{k-1,k,3}(n-2;(1))
-(\tilde{L}_{5}^{k-1,k,2}-\tilde{L}_{3}^{k-1,k,2})\no\\
&&+\frac{1}{2}\biggl(V_{1}^{k-1,k,1}(n-1;(1))\cdot
(\tilde{L}_{n-4}^{k-1,k,1}-\tilde{L}_{2}^{k-1,k,1})\no\\
&&+(\tilde{L}_{n-1}^{k-1,k,1}-\tilde{L}_{2}^{k-1,k,1})\cdot
V_{1}^{k-1,k,1}(n-3;(1))\no\\
&&-V_{1}^{k-1,k,4}(n-1;(1)+(2))\cdot
(\tilde{L}_{3}^{k-1,k,1}-\tilde{L}_{2}^{k-1,k,1})\no\\
&&-V_{1}^{k-1,k,4}(n-1;(3))\cdot
(\tilde{L}_{4}^{k-1,k,1}-\tilde{L}_{2}^{k-1,k,1})\biggr)
-V_{2}^{k-1,k,3}(6;(1))\no\\
&&+\frac{1}{4}(3hi_{1}(n)+3hi_{2}(n)+hi_{3}(n)),
\end{eqnarray}
\begin{eqnarray}
V_{3}^{k-1,k,5}(n;(1)+(1))&=&V_{3}^{k-1,k,4}(n;(1))+V_{3}^{k-1,k,4}(n-1;(1)) 
-(\tilde{L}_{6}^{k-1,k,3}-
\tilde{L}_{4}^{k-1,k,3})\no\\
&&+\frac{2}{3}\biggl(V_{2}^{k-1,k,3}(n;(1))\cdot(\tilde{L}_{n-4}^{k-1,k,1}-
\tilde{L}_{2}^{k-1,k,1})\no\\
&&+(\tilde{L}_{n}^{k-1,k,1}-
\tilde{L}_{2}^{k-1,k,1})\cdot V_{2}^{k-1,k,3}(n-2;(1))\no\\
&&-\bigl(V_{2}^{k-1,k,4}(n;(2))+V_{2}^{k-1,k,4}(n-1;(2))\no\\
&&-(\tilde{L}^{k-1,k,2}_{6}-
\tilde{L}^{k-1,k,2}_{3})\bigr)\cdot(\tilde{L}_{3}^{k-1,k,1}-
\tilde{L}_{2}^{k-1,k,1})\no\\
&&-V_{1}^{k-1,k,5}(n;(4))\cdot(\tilde{L}_{5}^{k-1,k,2}-
\tilde{L}_{3}^{k-1,k,2})\biggr)\no\\
&&+\frac{1}{3}\biggl(V_{1}^{k-1,k,2}(n;(1))\cdot(\tilde{L}_{n-3}^{k-1,k,2}-
\tilde{L}_{3}^{k-1,k,2})\no\\
&&+(\tilde{L}_{n}^{k-1,k,2}-
\tilde{L}_{3}^{k-1,k,2})\cdot V_{1}^{k-1,k,2}(n-3;(1))\no\\
&&-V_{1}^{k-1,k,5}(n;(1)+(3))\cdot
(\tilde{L}_{4}^{k-1,k,2}-\tilde{L}_{3}^{k-1,k,2})\no\\
&&-V_{2}^{k-1,k,5}(n;(3))
\cdot(\tilde{L}_{4}^{k-1,k,1}-\tilde{L}_{2}^{k-1,k,1})\biggr)\no\\
&&-\frac{1}{3}(\tilde{L}_{3}^{k-1,k,1}-\tilde{L}_{2}^{k-1,k,1})\cdot
(2hi_{1}(n)+hi_{2}(n)+hi_{3}(n)-hi_{4}(n))
\no\\     
\end{eqnarray}
The non-trivial rational coefficients of $V_{d-m}^{k-1,k,d}(n;\sigma_{m})$ come
from the rational factor, 
\begin{equation}
\prod_{j=1}^{l(\sigma_{m})-1}(1-\frac{n_{j}}{d-m}),\;\;(0\leq n_{j}<d-m),
\label{rat}
\end{equation} 
whose origin is the factor $\frac{1}{(d-m)^{l(\sigma_{m})-1}}$ in the 
definition of $V_{d-m}^{k-1,k,d}(n;\sigma_{m})$.
Then looking at the above formulas of $V_{5-m}^{k-1,k,5}(n;\sigma_{m})$
and $G_{5-m}^{k-1,k,5}(n;\sigma_{m})$, we can speculate that the rational 
factor in (\ref{rat}) is modified into, 
\begin{equation}
\prod_{j=1}^{l(\sigma_{m})-1}(1-\frac{n_{j}}{d}).
\label{rat2}
\end{equation}  
Let us take  $V_{2}^{k-1,k,5}(n;(1)+(1)+(1))$ as an example.
We decompose $V_{2}^{k-1,k,5}(n;(1)+(1)+(1))$ according to 
the decomposition in (\ref{dec}):
\begin{eqnarray}
V_{2}^{k-1,k,5}(n;(1)+(1)+(1))&=&\pi_{1}
\circ\pi_{1}\circ\pi_{1}
(\tilde{L}^{k-1,k,2}_{n}-\tilde{L}^{k-1,k,2}_{3})
+\pi_{1}\circ\pi_{1}(hi_2^{k-1,k,3}(n;(1)))\no\\
&&+\pi_{1}(hi_2^{k-1,k,4}(n;(1)+(1)))+hi_2^{k-1,k,5}(n;(1)+(1)+(1)).\no\\
\label{dec1}
\end{eqnarray}
More explicitly, each decomposed part corresponds to,
\begin{eqnarray}
&&\pi_{1}
\circ\pi_{1}\circ\pi_{1}
(\tilde{L}^{k-1,k,2}_{n}-\tilde{L}^{k-1,k,2}_{3})
+\pi_{1}\circ\pi_{1}(hi_2^{k-1,k,3}(n;(1)))=\no\\
&&V_{2}^{k-1,k,3}(n;(1))+V_{2}^{k-1,k,3}(n-1;(1))
-(\tilde{L}_{5}^{k-1,k,2}-\tilde{L}_{3}^{k-1,k,2})\no\\
&&+V_{2}^{k-1,k,3}(n-1;(1))+V_{2}^{k-1,k,3}(n-2;(1))
-(\tilde{L}_{5}^{k-1,k,2}-\tilde{L}_{3}^{k-1,k,2})
-V_{2}^{k-1,k,3}(6;(1)),\no\\
\end{eqnarray}
\begin{eqnarray}
&&\pi_{1}(hi_2^{k-1,k,4}(n;(1)+(1)))=\no\\
&&\frac{1}{2}\biggl(V_{1}^{k-1,k,1}(n;(1))\cdot
(\tilde{L}_{n-3}^{k-1,k,1}-\tilde{L}_{2}^{k-1,k,1})
+(\tilde{L}_{n}^{k-1,k,1}-\tilde{L}_{2}^{k-1,k,1})\cdot
V_{1}^{k-1,k,1}(n-2;(1))\no\\
&&-V_{1}^{k-1,k,4}(n;(1)+(2))\cdot
(\tilde{L}_{3}^{k-1,k,1}-\tilde{L}_{2}^{k-1,k,1})
-V_{1}^{k-1,k,4}(n;(3))\cdot
(\tilde{L}_{4}^{k-1,k,1}-\tilde{L}_{2}^{k-1,k,1})\no\\
&&+V_{1}^{k-1,k,1}(n-1;(1))\cdot
(\tilde{L}_{n-4}^{k-1,k,1}-\tilde{L}_{2}^{k-1,k,1})
+(\tilde{L}_{n-1}^{k-1,k,1}-\tilde{L}_{2}^{k-1,k,1})\cdot
V_{1}^{k-1,k,1}(n-3;(1))\no\\
&&-V_{1}^{k-1,k,4}(n-1;(1)+(2))\cdot
(\tilde{L}_{3}^{k-1,k,1}-\tilde{L}_{2}^{k-1,k,1})
-V_{1}^{k-1,k,4}(n-1;(3))\cdot
(\tilde{L}_{4}^{k-1,k,1}-\tilde{L}_{2}^{k-1,k,1})\biggr),\no\\
\end{eqnarray}
\begin{eqnarray}
&&hi_2^{k-1,k,5}(n;(1)+(1)+(1))
=\frac{1}{4}(3hi_{1}(n)+3hi_{2}(n)+hi_{3}(n)).
\label{exd}
\end{eqnarray}
According to (\ref{rat2}), the modification is given by,
\begin{equation}
\frac{1}{2}\rightarrow \frac{4}{5},\quad \frac{1}{4}=(\frac{1}{2})^{2}
\rightarrow (\frac{4}{5})^{2}=\frac{16}{25}.
\end{equation}
This modification is almost correct, but there exist some errors in
the modification of (\ref{exd}), which are given by, 
\begin{equation}
-\frac{2}{25}(hi_{1}(n)+hi_{2}(n)+hi_{4}(n)).
\end{equation}
Similar errors also occur in the cases of $V_{2}^{k-1,k,5}(n;(1)+(2))$
and $V_{3}^{k-1,k,5}(n;(1)+(1))$. One of the reasons of such errors
comes from the fact
that the decomposition of (\ref{dec}) is not unique, as was suggested in 
the remark of (\ref{dec}), but there must be other reasons because in
the case of $V_{2}^{k-1,k,5}(n;(1)+(1)+(1))$, the decomposition is
unique. Therefore, further consideration is needed.
\begin{Q} 
Fix the general rule of the modification of 
$V_{d-m}^{N,k,d}(n;\sigma_{m})$ into $G_{d-m}^{N,k,d}(n;\sigma_{m})$.
\end{Q} 
\begin{Rem}
$\pi_{1}(hi_{2}^{k-1,k,4}(n;(1)+(1)))$ does not vanish when $n=7$, but 
Conjecture 4 tells us that it receives modification. Thus, the condition
(iv) of Proposition 1 does not hold true for $G_{d-m}^{N,k,d}(n;\sigma_{m})$
in the $d\geq 5$ cases.
\end{Rem}

{\bf Acknowledgement}\\
We would like to thank T.Eguchi, S.Hosono and A.Collino for discussions.
We also thank M.Naka for kind encouragement. Research of the author is 
supported by JSPS postodoctoral fellowship.


\newpage 
\begin{thebibliography}{99}
\bibitem{beauville}A.Beauville.
\newblock{\it Quantum Cohomology of Complete Intersections}
\newblock Mathematical, Physics Analysis and Geometry  168 (1995),
 384-398.
\bibitem{bert}A.Bertram.
\newblock{\em Another way to enumerate rational curves with torus
 actions}
\newblock{Preprint, alg-geom/9905159}
\bibitem{harris} L.Caporaso, J.Harris. 
\newblock{\em Counting Plane Curves of Any Genus}
\newblock alg-geom/9608025.
\bibitem{cj} A.Collino, M.Jinzenji.
\newblock{\em On the Structure of the Small Quantum Cohomology 
Rings of Projective Hypersurfaces}
\newblock{Commun. Math. Phys. 206 (1999) 157.}
\bibitem{gath} A. Gathmann.
{\em  Absolute and Relative Gromov-Witten Invariants of Very Ample 
Hypersurfaces}
\newblock math.AG/9908054.
\bibitem{giv}A.B.Givental.
\newblock{\em Equivariant Gromov-Witten Invariants}
\newblock Internat. Math. Res.Notices 13 (1996),613--663.
\bibitem{gmp}B.R.Greene, D.R.Morrison and M.R.Plesser.
\newblock{\em Mirror Manifolds in Higher Dimension}
\newblock{Commun.Math.Phys. 173 (1995) 559-598.} 
\bibitem{6}M.Jinzenji.
\newblock{\em Completion of the Conjecture: Quantum Cohomology Rings of 
Fano Hypersurfaces}
\newblock To appear in Mod. Phys. Lett. A.
\bibitem{jin}M.Jinzenji.
\newblock{\em On Quantum Cohomology Rings for Hypersurfaces in
 $CP^{N-1}$}
\newblock {J. Math. Phys. 38 6613-6638.1997.}
\bibitem{gene}M.Jinzenji.
\newblock{\em On the Quantum Cohomology Rings of General Type Projective 
Hypersurfaces and Generalized Mirror Transformation}
\newblock{To appear in Int. J. of Mod. Phys. A.}
\bibitem{torus} M.Kontsevich.
\newblock{\em Enumeration of Rational Curves via Torus Actions}
\newblock{The moduli space of curves, R.Dijkgraaf, C.Faber, G.van der
 Geer (Eds.), Progress in Math., v.129, Birkh\"auser, 1995, 335-368.}
\bibitem{km}M.Kontsevich, Y.Manin.
\newblock{\em Gromov-Witten Classes, Quantum Cohomology, and Enumerative 
Geometry.}
\newblock{Commun. Math. Phys. 164, 525-562 (1994)}
\bibitem{yau}B.Lian, K.Liu and S.T.Yau.
\newblock{\em Mirror Principle III}
\newblock{math-AG/9912038.}
\bibitem{rt}Y.Ruan, G.Tian.
\newblock{\em A Mathematical Theory of Quantum Cohomology.}
\newblock{J. Diff. Geom. 42 no.2, (1995)}
\end{thebibliography}
\end{document}